\documentclass[12pt]{article}
\usepackage{amssymb}
\usepackage{amsfonts}
\usepackage{amsmath}

\input{tcilatex}
\begin{document}

\begin{center}
\textbf{DIAMETER\ DIMINISHING\ TO\ ZERO\ IFSs}

\bigskip

by Radu MICULESCU\ and Alexandru MIHAIL

\bigskip
\end{center}

\textbf{Abstract}. {\small In this paper we introduce the notion of diameter
diminishing to zero iterated function system, study its properties and
provide alternative characterizations of it.}

\bigskip



\bigskip

\textbf{1. INTRODUCTION}

\bigskip

The notion of iterated function system, which is due to J. Hutchinson (see
[9]), was popularized by M. Barnsley (see [3]). It represents one of the
main ways to generate fractal sets and (since it has numerous applications)
various generalizations of this concept were introduced. Among them we
mention the one exhibited by A. Kameyama (see [11]) under the label of
self-similar system which is a topological generalization of the attractor
of an iterated function system.

In connection with Kamayema's work, R. Atkins, M. Barnsley, A. Vince and D.
Wilson (see [1]) presented a theorem that characterizes hyperbolic affine
iterated function systems defined on $\mathbb{R}^{m}$. Results connected
with the aforementioned theorem are included in [2], [4], [13], [14], [15],
[16], [18], [22] and [23].

Along these lines of research, in this paper, we introduce the concept of
(hyperbolic) diameter diminishing to zero iterated function system (see
Definition 2.17 and Remarks 2.19) and study its properties (see Propositions
3.1 - 3.7).

In addition, via the concepts of hyperbolic $\varphi $-contractive iterated
function system (see Definition 2.11 and Remarks 2.19), hyperbolic (locally)
uniformly point fibred iterated function system (see Definitions 2.13 and
2.14 and Remarks 2.19) and iterated function system having hyperbolic
attractor (see Definition 2.15 and Remarks 2.19), we provide alternative
characterizations of hyperbolic diameter diminishing to zero iterated
function systems (see Theorem 4.1). Consequently we also come up with
different ways to prove that an iterated function system has hyperbolic
attractor.

\bigskip

\textbf{2.} \textbf{PRELIMINARIES}

\bigskip

Given two sets $A$ and $B$, by $B^{A}$ we mean the set of functions from $A$
to $B$.

By $\mathbb{N}$ we mean the set $\{0,1,2,...,n,...\}$ and by $\mathbb{N}%
^{\ast }$ we mean the set $\{1,2,...,n,...\}$.

For a set $X$, a function $f:X\rightarrow X$ and $n\in \mathbb{N}^{\ast }$,
by $f^{[n]}$ we mean the composition of $f$ by itself $n$ times. By $f^{[0]}$
we mean the identity function $Id_{X}:X\rightarrow X$ given by $Id_{X}(x)=x$
for every $x\in X$.

Given a metric space $(X,d)$, by:

- $P_{b}(X)$ we mean the set of non-empty bounded subsets of $X$

- $P_{b,cl}(X)$ we mean the set of non-empty bounded and closed subsets of $%
X $

- $P_{cp}(X)$ we mean the set of non-empty compact subsets of $X$

- by $B(x,r)$ we mean the set $\{y\in X\mid d(x,y)<r\}$, where $x\in X$ and $%
r>0$.

\bigskip

\textbf{The Hausdorff-Pompeiu metric}

\bigskip

\textbf{Definition 2.1.} \textit{Given a metric space }$(X,d)$\textit{, }$%
H_{d}:P_{b,cl}(X)\times P_{b,cl}(X)\rightarrow \lbrack 0,+\infty )$\textit{\
given by}

\begin{equation*}
H_{d}(A,B)=\max \{\underset{x\in A}{\sup }\text{ }d(x,B),\underset{x\in B}{%
\sup }\text{ }d(x,A)\}=
\end{equation*}%
\begin{equation*}
=\inf \{\varepsilon \in \lbrack 0,\infty )\mid A\subseteq E_{\varepsilon }(B)%
\text{ \textit{and} }B\subseteq E_{\varepsilon }(A)\}\text{,}
\end{equation*}%
\textit{for all }$A,B\in P_{b,cl}(X)$\textit{, where }%
\begin{equation*}
d(x,A)=\underset{y\in A}{\inf }\text{ }d(x,y)
\end{equation*}%
\textit{and} 
\begin{equation*}
E_{\varepsilon }(A)\overset{def}{=}\{y\in X\mid \text{\textit{there exists}}%
\ x\in A\ \text{\textit{such that\ }}d(x,y)<\varepsilon \}=\underset{x\in A}{%
\cup }B(x,\varepsilon )\text{\textit{,}}
\end{equation*}%
\textit{\ turns out to be a metric, which is called the Hausdorff-Pompeiu
metric.}

\bigskip

\textbf{Proposition 2.2.} \textit{For a metric space }$(X,d)$\textit{, we
have}

\begin{equation*}
H_{d}(\underset{i\in I}{\cup }H_{i},\underset{i\in I}{\cup }K_{i})\leq 
\underset{i\in I}{\sup }H_{d}(H_{i},K_{i})\text{,}
\end{equation*}%
\textit{for every} $(H_{i})_{i\in I}$\textit{\ and }$(K_{i})_{i\in I}$%
\textit{\ families of elements from }$P_{b,cl}(X)$ \textit{such that} $%
\underset{i\in I}{\cup }H_{i},\underset{i\in I}{\cup }K_{i}\in P_{b,cl}(X)$.

\bigskip

\textbf{Remark 2.3. }\textit{The metric spaces }$(P_{b,cl}(X),H_{d})$\textit{%
\ and }$(P_{cp}(X),H_{d})$\textit{\ are complete provided that }$(X,d)$%
\textit{\ is complete.}

\bigskip

\textbf{The shift space}

\bigskip

Let $I$ be a non-empty set.

We denote the set $I^{\mathbb{N}^{\ast }}$ by $\Lambda (I)$. Thus $\Lambda
(I)$ is the set of infinite words with letters from the alphabet $I$ and a
standard element $\omega $ of $\Lambda (I)$ can be presented as $\omega
=\omega _{1}\omega _{2}...\omega _{n}\omega _{n+1}...$ .

We denote the set $I^{\{1,2,...,n\}}$ by $\Lambda _{n}(I)$. Thus $\Lambda
_{n}(I)$ is the set of words with letters from the alphabet $I$ of length $n$
and a standard element $\omega $ of $\Lambda _{n}(I)$ can be presented as $%
\omega =\omega _{1}\omega _{2}...\omega _{n}$. The length of $\omega $ is
denoted by $\left\vert \omega \right\vert $. By $\Lambda _{0}(I)$ we mean
the set having only one element, namely the empty word denoted by $\lambda $.

We denote the set $\underset{n\in \mathbb{N}}{\cup }\Lambda _{n}(I)$ by $%
\Lambda ^{\ast }(I)$. Thus $\Lambda ^{\ast }(I)$ is the set of words with
letters from the alphabet having finite length.

Given $m,n\in \mathbb{N}$ and two words $\omega =\omega _{1}\omega
_{2}...\omega _{n}\in \Lambda _{n}(I)$\ and $\theta =\theta _{1}\theta
_{2}...\theta _{m}\in \Lambda _{m}(I)$ or $\theta =\theta _{1}\theta
_{2}...\theta _{m}\theta _{m+1}...\in \Lambda (I)$, by $\omega \theta $ we
mean the concatenation of the words $\omega $ and $\theta $, i.e.$\ \omega
\theta =\omega _{1}\omega _{2}...\omega _{n}\theta _{1}\theta _{2}...\theta
_{m}$ and respectively $\omega \theta =\omega _{1}\omega _{2}...\omega
_{n}\theta _{1}\theta _{2}...\theta _{m}\theta _{m+1}...$ .

For $m\in \mathbb{N}^{\ast }$ and $\omega =\omega _{1}\omega _{2}...\omega
_{n}\omega _{n+1}...\in \Lambda (I)$, by $[\omega ]_{m}$ we mean $\omega
_{1}\omega _{2}...\omega _{m}$.

For\textit{\ }$\omega =\omega _{1}\omega _{2}...\omega _{n}\in \Lambda
_{n}(I)$, the word $\theta =\theta _{1}\theta _{2}...\theta _{m}\theta
_{m+1}...\in \Lambda (I)$, where $\theta _{nk+i}=\omega _{i}$ for every $%
k\in \mathbb{N}$ and every $i\in \{1,...,n-1,n\}$, will be denoted by $%
\overset{.}{\omega }$.

For $i\in I$ one can consider the function $\tau _{i}:\Lambda (I)\rightarrow
\Lambda (I)$ given by%
\begin{equation*}
\tau _{i}(\omega )=i\omega \text{,}
\end{equation*}%
for every $\omega \in \Lambda (I)$.

$\Lambda (I)$ becomes a metric space if we endow it with the distance
described by $d_{\Lambda }(\omega ,\theta )=\{%
\begin{array}{cc}
0\text{,} & \text{if }\omega =\theta \\ 
\frac{1}{2^{\min \{k\in \mathbb{N}^{\ast }\mid \omega _{k}\neq \theta _{k}\}}%
}\text{,} & \text{if }\omega \neq \theta%
\end{array}%
$, where $\omega =\omega _{1}\omega _{2}\omega _{3}...\omega _{n}\omega
_{n+1}...$ and $\theta =\theta _{1}\theta _{2}\theta _{3}...\theta
_{n}\theta _{n+1}...$.

\newpage

\textbf{Remarks 2.4.}

a)\textit{\ The convergence in the metric space }$(\Lambda (I),d_{\Lambda })$%
\textit{\ is the convergence on components.}

b)\textit{\ }$(\Lambda (I),d_{\Lambda })$\textit{\ is a complete metric
space.}

c)\textit{\ If }$I$\textit{\ is finite, then} $(\Lambda (I),d_{\Lambda })$%
\textit{\ is compact.}

\bigskip

Note that points a) and c) follow from the fact that the metric $d_{\Lambda
} $ induces the Tychonoff product topology.

\bigskip

\textbf{Proposition 2.5.} \textit{Let us consider} $(\omega _{n})_{n\in 
\mathbb{N}^{\substack{ \ast  \\  \\  \\ }}}\subseteq \Lambda (I)$ \textit{and%
} $\omega \in \Lambda (I)$ \textit{such that }$\underset{n\rightarrow \infty 
}{\lim }\omega _{n}=\omega $\textit{. Then for every} $m\in \mathbb{N}^{\ast
}$ \textit{there exists} $n_{m}\in \mathbb{N}^{\ast }$ \textit{such that }$%
[\omega _{n}]_{m}=[\omega ]_{m}$\textit{\ for every} $n\in \mathbb{N}^{\ast }
$, $n\geq n_{m}$\textit{.}

\textit{Proof}. Ignoring those $\omega _{n}$ which are equal with $\omega $
and supposing that $\omega _{n}=\omega _{n}^{1}\omega _{n}^{2}...\omega
_{n}^{k}...$ and $\omega =\omega _{1}\omega _{2}...\omega _{k}...$, for
every $m\in \mathbb{N}^{\ast }$ there exists $n_{m}\in \mathbb{N}^{\ast }$
such that 
\begin{equation*}
\frac{1}{2^{\min \{l\in \mathbb{N}^{\ast }\mid \omega _{n}^{l}\neq \omega
_{l}\}}}<\frac{1}{2^{m}}\text{,}
\end{equation*}%
for every $n\in \mathbb{N}^{\ast }$, $n\geq n_{m}$. As the previous
inequality means that 
\begin{equation*}
m<\min \{l\in \mathbb{N}^{\ast }\mid \omega _{n}^{l}\neq \omega _{l}\}\text{,%
}
\end{equation*}%
we get $\omega _{n}^{1}=\omega _{1},\omega _{n}^{2}=\omega _{2},...,\omega
_{n}^{m}=\omega _{m}$, i.e 
\begin{equation*}
\lbrack \omega _{n}]_{m}=[\omega ]_{m}\text{,}
\end{equation*}%
for every $n\in \mathbb{N}^{\ast }$, $n\geq n_{m}$. $\square $

\bigskip

\textbf{Comparison functions and }$\varphi $\textbf{-contractions}

\bigskip

\textbf{Definition 2.6 }(comparison function)\textbf{.} \textit{A function }$%
\varphi :[0,\infty )\rightarrow \lbrack 0,\infty )$\textit{\ is called a
comparison function if:}

\textit{i) }$\varphi $\textit{\ is increasing;}

\textit{ii) }$\varphi (t)<t$\textit{\ for every }$t>0$\textit{;}

\textit{iii) }$\varphi $\textit{\ is right-continuous.}

\bigskip

\textbf{Remark 2.7.} \textit{If }$\varphi :[0,\infty )\rightarrow \lbrack
0,\infty )$\textit{\ is a comparison function, then }$\underset{n\rightarrow
\infty }{\lim }\varphi ^{\lbrack n]}(t)=0$\textit{\ for every }$t\in \lbrack
0,\infty )$\textit{.}

\bigskip

\textbf{Definition 2.8 }($\varphi $-contraction)\textbf{.} \textit{Given a
metric space }$(X,d)$\textit{\ and a comparison function }$\varphi $\textit{%
, a function }$f:X\rightarrow X$ \textit{is called }$\varphi $\textit{%
-contraction if }$d(f(x),f(y))\leq \varphi (d(x,y))$\textit{\ for all }$%
x,y\in X$\textit{.}

\bigskip

A very good reference on comparison functions and $\varphi $-contractions is
the survey [10]. See also [8].

\bigskip

\textbf{Iterated function systems}

\bigskip

\textbf{Definition 2.9 }(iterated function system)\textbf{.} \textit{A pair} 
$((X,d),(f_{i})_{i\in I})$ \textit{is called an iterated function system }%
(IFS for short) \textit{if:}

\textit{i) }$(X,d)$\textit{\ is a complete metric space;}

\textit{ii) }$I$\textit{\ is a finite set;}

\textit{iii)\ }$f_{i}:X\rightarrow X$\textit{\ is continuous for each }$i\in
I$\textit{;}

\textit{iv)} $f_{i}(B)\in P_{b}(X)$ \textit{for every }$B\in P_{b}(X)$%
\textit{\ and every} $i\in I$.

\bigskip

\textbf{Notations}

\textbf{1}. \textit{We shall denote the IFS }$((X,d),(f_{i})_{i\in I})$%
\textit{\ by} $\mathcal{S}$.

\textbf{2}. \textit{In the framework of the above definition, for }$\omega
=\omega _{1}\omega _{2}...\omega _{n}\in \Lambda _{n}(I)$\textit{\ and }$B$%
\textit{\ subset of }$X$\textit{, by }$f_{\omega }(B)$\textit{\ we mean }$%
(f_{\omega _{1}}\circ ...\circ f_{\omega _{n}})(B)$\textit{.}

\bigskip

\textbf{Definition 2.10} (fractal operator)\textbf{.} \textit{Given an IFS} $%
\mathcal{S=}((X,d),(f_{i})_{i\in I})$\textit{, the function }$F_{\mathcal{S}%
}:P_{b,cl}(X)\rightarrow P_{b,cl}(X)$\textit{, given by}%
\begin{equation*}
F_{\mathcal{S}}(B)=\overline{\underset{i\in I}{\cup }f_{i}(B)}\text{,}
\end{equation*}%
\textit{\ for every }$B\in P_{b,cl}(X)$\textit{, is called the fractal
operator associated to} $\mathcal{S}$.

\bigskip

\textbf{Definition 2.11 }($\varphi $-contractive IFS)\textbf{.} \textit{%
Given a comparison function} $\varphi $, \textit{an iterated function system 
}$\mathcal{S}=((X,d),(f_{i})_{i\in I})$\textit{\ is called }$\varphi $%
\textit{-contractive if }$f_{i}$\textit{\ is a }$\varphi $\textit{%
-contraction for each }$i\in I$\textit{.}

\bigskip

\textbf{Definition 2.12 }(point fibred IFS)\textbf{.} \textit{An iterated
function system }$\mathcal{S}=((X,d),(f_{i})_{i\in I})$\textit{\ is called
point fibred if for every} $\omega \in \Lambda (I)$ \textit{there exists }$%
a_{\omega }\in X$\textit{\ such that}%
\begin{equation*}
\underset{n\rightarrow \infty }{\lim }f_{[\omega ]_{n}}(x)=a_{\omega }\text{,%
}
\end{equation*}%
\textit{for all }$x\in X$\textit{.}

\bigskip

\textbf{Definition 2.13 }(locally uniformly point fibred IFS)\textbf{.} 
\textit{An iterated function system }$\mathcal{S}=((X,d),(f_{i})_{i\in I})$%
\textit{\ is called locally uniformly point fibred if it is point fibred and
for each }$x\in X$ \textit{there exists an open set }$D_{x}$\textit{\
containing }$x$\textit{\ such that }

\begin{equation*}
\underset{n\rightarrow \infty }{\lim }\underset{\omega \in \Lambda (I)}{\sup 
}\underset{y\in D_{x}}{\sup }d(f_{[\omega ]_{n}}(y),a_{\omega })=0\text{.}
\end{equation*}

\bigskip

Note that the concept of locally uniformly point fibred IFS is the same as
the "condition C" that was introduced in Definition 3.1 from [15].

\bigskip

\textbf{Definition 2.14 }(uniformly point fibred IFS)\textbf{.} \textit{An
iterated function system }$\mathcal{S}=((X,d),(f_{i})_{i\in I})$\textit{\ is
called uniformly point fibred if it is point fibred and }%
\begin{equation*}
\underset{n\rightarrow \infty }{\lim }\underset{\omega \in \Lambda (I)}{\sup 
}\underset{x\in B}{\sup }d(f_{[\omega ]_{n}}(x),a_{\omega })=0\text{,}
\end{equation*}%
\textit{for every} $B\in P_{b,cl}(X)$.

\bigskip

\textbf{Definition 2.15 }(IFS having attractor)\textbf{.} \textit{We say
that an iterated function system }$\mathcal{S}=((X,d),(f_{i})_{i\in I})$%
\textit{\ has an attractor if there exists }$A_{\mathcal{S}}\in P_{b,cl}(X)$%
\textit{\ such that:}

\textit{i)} 
\begin{equation*}
F_{\mathcal{S}}(A_{\mathcal{S}})=A_{\mathcal{S}}\text{;}
\end{equation*}

\textit{ii)} 
\begin{equation*}
\underset{n\rightarrow \infty }{\lim }H_{d}(F_{\mathcal{S}}^{[n]}(B),A_{%
\mathcal{S}})=0\text{,}
\end{equation*}
\textit{for each} $B\in P_{b,cl}(X)$.

\bigskip

\textbf{Remark 2.16.} $F_{\mathcal{S}}$ \textit{has a unique fixed point,
namely} $A_{\mathcal{S}}\in P_{cp}(X)$\textit{, which is called the
attractor of }$\mathcal{S}$.

Indeed, if for some $B\in P_{b,cl}(X)$ we have $F_{\mathcal{S}}(B)=B$, then 
\begin{equation*}
0\overset{\text{ii)}}{=}\underset{n\rightarrow \infty }{\lim }H_{d}(F_{%
\mathcal{S}}^{[n]}(B),A_{\mathcal{S}})=\underset{n\rightarrow \infty }{\lim }%
H_{d}(B,A_{\mathcal{S}})=H_{d}(B,A_{\mathcal{S}})\text{,}
\end{equation*}%
so $B=A_{\mathcal{S}}$. In addition, $\underset{n\rightarrow \infty }{\lim }%
H_{d}(F_{\mathcal{S}}^{[n]}(K),A_{\mathcal{S}})=0$ for every $K\in
P_{cp}(X)\subseteq P_{b,cl}(X)$ and in view of Proposition 2.7 from [17] we
conclude that $A_{\mathcal{S}}\in P_{cp}(X)$.

\bigskip

\textbf{Definition 2.17 }(diameter diminishing to zero IFS)\textbf{.} 
\textit{An iterated function system }$\mathcal{S}=((X,d),(f_{i})_{i\in I})$%
\textit{\ is called diameter diminishing to zero iterated function systems
if for every }$B\in P_{b,cl}(X)$ \textit{there exists} $M_{B}\in P_{b,cl}(X)$
\textit{such that:}

\textit{i)} 
\begin{equation*}
B\subseteq M_{B}\text{;}
\end{equation*}

\textit{ii)}%
\begin{equation*}
F_{\mathcal{S}}(M_{B})\subseteq M_{B}\text{;}
\end{equation*}

\textit{iii)}%
\begin{equation*}
\underset{n\rightarrow \infty }{\lim }\text{ }\underset{\omega \in \Lambda
_{n}(I)}{\max }diam(f_{\omega }(M_{B}))=0\text{.}
\end{equation*}

\bigskip

\textbf{Definition 2.18 }(hyperbolic $\varphi $-contractive IFS)\textbf{.} 
\textit{Given a comparison function} $\varphi $, \textit{an iterated
function system }$\mathcal{S}=((X,d),(f_{i})_{i\in I})$\textit{\ is called
hyperbolic }$\varphi $\textit{-contractive if there exists a distance }$%
d_{1} $ \textit{on }$X$\textit{\ such that:}

\textit{i) }$d$\textit{\ and }$d_{1}$\textit{\ are topologically equivalent;}

\textit{ii) }$(X,d_{1})$\textit{\ is complete;}

\textit{iii)} $\mathcal{S}_{1}=((X,d_{1}),(f_{i})_{i\in I})$ \textit{is }$%
\varphi $\textit{-contractive.}

\bigskip

\textbf{Remarks 2.19.}

a) \textit{The concepts of hyperbolic locally uniformly point fibred IFS,
hyperbolic uniformly point fibred IFS, IFS having hyperbolic attractor and
hyperbolic diameter diminishing to zero IFS could be defined having as model
Definition 2.18.}

b) \textit{An iterated function system }$\mathcal{S}$ \textit{which is
uniformly point fibred is locally uniformly point fibred.}

c) \textit{An iterated function system }$\mathcal{S}$ \textit{which is
hyperbolic uniformly point fibred is hyperbolic locally uniformly point
fibred.}

d) \textit{As one of the referees of this paper noted, Definition 2.15
raises the question whether there exists and IFS with hyperbolic attractor
and not having attractor.}

\bigskip

\textbf{3. THE\ PROPERTIES\ OF\ DIAMETER\ DIMINISHING\ TO\ ZERO\ ITERATED\
FUNCTION\ SYSTEMS}

\bigskip

\textbf{Proposition 3.1.} \textit{Given a} \textit{diameter diminishing to
zero IFS} $\mathcal{S}=((X,d),(f_{i})_{i\in I})$, \textit{for every }$\omega
\in \Lambda (I)$\textit{\ there exists }$a_{\omega }\in X$\textit{\ such that%
}%
\begin{equation*}
\underset{n\in \mathbb{N}^{\ast }}{\cap }\overline{f_{[\omega ]_{n}}(B)}%
=\{a_{\omega }\}\text{,}
\end{equation*}%
\textit{for every }$B\in P_{b,cl}(X)$\textit{\ having the property that }$F_{%
\mathcal{S}}(B)\subseteq B$\textit{.}

\textit{Proof}. Let $\omega \in \Lambda (I)$ be fixed.

\textbf{Claim 1}. For every $B\in P_{b,cl}(X)$ such that $F_{\mathcal{S}%
}(B)\subseteq B$ there exists $a_{\omega ,B}\in X$ such that 
\begin{equation*}
\underset{n\in \mathbb{N}^{\ast }}{\cap }\overline{f_{[\omega ]_{n}}(B)}%
=\{a_{\omega ,B}\}\text{.}
\end{equation*}

\textit{Justification of Claim 1}. As $f_{i}(B)\subseteq B$ for each $i\in I$%
, we have 
\begin{equation*}
f_{[\omega ]_{n+1}}(B)\subseteq f_{[\omega ]_{n}}(B)\text{,}
\end{equation*}%
for every $n\in \mathbb{N}^{\ast }$. Since%
\begin{equation*}
diam(\overline{f_{[\omega ]_{n}}(B)})=diam(f_{[\omega ]_{n}}(B))\leq 
\end{equation*}%
\begin{equation}
\leq \underset{\theta \in \Lambda _{n}(I)}{\max }diam(f_{\theta }(B))\overset%
{B\subseteq M_{B}}{\leq }\underset{\theta \in \Lambda _{n}(I)}{\max }%
diam(f_{\theta }(M_{B}))\text{,}  \tag{1}
\end{equation}%
for every $n\in \mathbb{N}^{\ast }$, where $M_{B}$ is the element of $%
P_{b,cl}(X)$ whose existence is stated in Definition 2.17. Using iii) of the
same Definition, via $(1)$, we conclude that%
\begin{equation*}
\underset{n\rightarrow \infty }{\lim }diam(\overline{f_{[\omega ]_{n}}(B)})=0%
\text{.}
\end{equation*}%
Therefore, according to Cantor's theorem, there exists $a_{\omega ,B}\in X$
such that%
\begin{equation*}
\underset{n\in \mathbb{N}^{\ast }}{\cap }\overline{f_{[\omega ]_{n}}(B)}%
=\{a_{\omega ,B}\}\text{.}
\end{equation*}%
The justification of the claim is done.

\textbf{Claim 2}. For every $B,C\in P_{b,cl}(X)$ such that $F_{\mathcal{S}%
}(B)\subseteq B$, $F_{\mathcal{S}}(C)\subseteq C$ and $B\subseteq C$, we
have 
\begin{equation*}
a_{\omega ,B}=a_{\omega ,C}\text{.}
\end{equation*}

\textit{Justification of Claim 2}. Indeed 
\begin{equation*}
\{a_{\omega ,B}\}\overset{\text{Claim 1}}{=}\underset{n\in \mathbb{N}^{\ast }%
}{\cap }\overline{f_{[\omega ]_{n}}(B)}\subseteq \underset{n\in \mathbb{N}%
^{\ast }}{\cap }\overline{f_{[\omega ]_{n}}(C)}\overset{\text{Claim 1}}{=}%
\{a_{\omega ,C}\}
\end{equation*}%
and the justification of the claim is done.

\textbf{Claim 3}. For every $B,C\in P_{b,cl}(X)$ such that $F_{\mathcal{S}%
}(B)\subseteq B$ and $F_{\mathcal{S}}(C)\subseteq C$, we have 
\begin{equation*}
a_{\omega ,B}=a_{\omega ,C}\text{.}
\end{equation*}

\textit{Justification of Claim 3}. Since 
\begin{equation*}
B\cup C\in P_{cl,b}(X)
\end{equation*}%
and 
\begin{equation*}
F_{\mathcal{S}}(B\cup C)\subseteq F_{\mathcal{S}}(B)\cup F_{\mathcal{S}%
}(C)\subseteq B\cup C\text{,}
\end{equation*}%
we infer that 
\begin{equation*}
a_{\omega ,B}\overset{\text{Claim 2}}{=}a_{\omega ,B\cup C}\overset{\text{%
Claim 2}}{=}a_{\omega ,C}\text{,}
\end{equation*}%
and the justification of the claim is done.

Hence $\{a_{\omega ,B}\mid B\in P_{b,cl}(X)$ such that $F_{\mathcal{S}%
}(B)\subseteq B\}$ has only one element, which is denoted by $a_{\omega }$. $%
\square $

\bigskip

In view of the previous Proposition, given a diameter diminishing to zero
IFS $\mathcal{S}=((X,d),(f_{i})_{i\in I})$ one can consider the function $%
\pi :\Lambda (I)\rightarrow X$ given by 
\begin{equation*}
\pi (\omega )=a_{\omega }\text{,}
\end{equation*}%
\ for each $\omega \in \Lambda (I)$.

\bigskip

\textbf{Proposition 3.2.} \textit{Each diameter diminishing to zero IFS} 
\textit{is} \textit{uniformly point fibred.}

\textit{Proof}. If $\mathcal{S}=((X,d),(f_{i})_{i\in I})$, then, according
to Definition 2.17, for each $B\in P_{b,cl}(X)$ exists $M_{B}\in P_{b,cl}(X)$
such that $B\subseteq M_{B}$ and $F_{\mathcal{S}}(M_{B})\subseteq M_{B}$.
Consequently, for every $\omega \in \Lambda (I)$, we have%
\begin{equation*}
\underset{x\in B}{\sup }\text{ }d(f_{[\omega ]_{n}}(x),a_{\omega })\leq 
\underset{x\in M_{B}}{\sup }d(f_{[\omega ]_{n}}(x),a_{\omega })\overset{%
a_{\omega }\in \overline{f_{[\omega ]_{n}}(M_{B})}}{\leq }
\end{equation*}%
\begin{equation*}
\leq diam(\overline{f_{[\omega ]_{n}}(M_{B})})\leq \underset{\theta \in
\Lambda _{n}(I)}{\max }diam(f_{\theta }(M_{B}))\text{,}
\end{equation*}%
for every $n\in \mathbb{N}^{\ast }$. Hence%
\begin{equation*}
\underset{\omega \in \Lambda (I)}{\sup }\underset{x\in B}{\sup }\text{ }%
d(f_{[\omega ]_{n}}(x),a_{\omega })\leq \underset{\omega \in \Lambda _{n}(I)}%
{\max }diam(f_{\omega }(M_{B}))\text{,}
\end{equation*}%
for every $n\in \mathbb{N}^{\ast }$, so, in view of Definition 2.17, we
obtain%
\begin{equation*}
\underset{n\rightarrow \infty }{\lim }\underset{\omega \in \Lambda (I)}{\sup 
}\underset{x\in B}{\sup }d(f_{[\omega ]_{n}}(x),a_{\omega })=0\text{,}
\end{equation*}%
i.e. $\mathcal{S}$ is uniformly point fibred. $\square $

\bigskip

\textbf{Proposition 3.3.} \textit{For} \textit{each diameter diminishing to
zero IFS }$\mathcal{S}=((X,d),(f_{i})_{i\in I})$\textit{, we have} 
\begin{equation*}
f_{i}\circ \pi =\pi \circ \tau _{i}\text{,}
\end{equation*}%
\textit{for every }$i\in I$.

\textit{Proof}. For a fixed $B\in P_{b,cl}(X)$ such that $F_{\mathcal{S}%
}(B)\subseteq B$, we have%
\begin{equation*}
\{f_{i}(\pi (\omega ))\}=\{f_{i}(a_{\omega })\}=f_{i}(\underset{n\in \mathbb{%
N}^{\ast }}{\cap }\overline{f_{[\omega ]_{n}}(B)})\overset{f_{i}\text{
continuous}}{\subseteq }\underset{n\in \mathbb{N}^{\ast }}{\cap }\overline{%
f_{i}(f_{[\omega ]_{n}}(B))}=
\end{equation*}%
\begin{equation*}
=\underset{n\in \mathbb{N}^{\ast }}{\cap }\overline{f_{[i\omega ]_{n}}(B)}%
=\{a_{i\omega }\}=\{\pi (i\omega )\}=\{\pi (\tau _{i}(\omega ))\}\text{,}
\end{equation*}%
i.e.%
\begin{equation*}
(f_{i}\circ \pi )(\omega )=(\pi \circ \tau _{i})(\omega )\text{,}
\end{equation*}%
for every $\omega \in \Lambda (I)$ and every $i\in I$. $\square $

\bigskip

\textbf{Proposition 3.4.} \textit{For} \textit{each diameter diminishing to
zero IFS }$\mathcal{S}=((X,d),(f_{i})_{i\in I})$\textit{, the function} $\pi 
$\textit{\ is continuous.}

\textit{Proof.} First of all, let us chose a fixed $B\in P_{b,cl}(X)$ such
that $F_{\mathcal{S}}(B)\subseteq B$. Now let us consider $\omega \in
\Lambda (I)$ and a sequence $(\omega _{n})_{n\in \mathbb{N}^{\ast }}$ of
elements of $\Lambda (I)$ converging to $\omega $. Therefore, in view of
Proposition 2.5, for each $m\in \mathbb{N}^{\ast }$ there exists $n_{m}\in 
\mathbb{N}^{\ast }$ such that 
\begin{equation}
\lbrack \omega _{n}]_{m}=[\omega ]_{m}\text{,}  \tag{2}
\end{equation}%
for every $n\in \mathbb{N}^{\ast }$, $n\geq n_{m}$. Let us consider a fixed,
but arbitrarily chosen $\varepsilon >0$. Taking into account Definition
2.17, there exists $m_{\varepsilon }\in \mathbb{N}^{\ast }$ such that 
\begin{equation}
\underset{\theta \in \Lambda _{m_{\varepsilon }}(I)}{\max }diam(f_{\theta
}(B))<\varepsilon \text{.}  \tag{3}
\end{equation}%
Then for every $n\in \mathbb{N}^{\ast }$, $n\geq n_{m_{\varepsilon }}$, as $%
\pi (\omega _{n})\in \overline{f_{[\omega _{n}]_{m_{\varepsilon }}}(B)}$ and 
$\pi (\omega )\in \overline{f_{[\omega ]_{m_{\varepsilon }}}(B)}$, with the
notation $[\omega _{n}]_{m_{\varepsilon }}\overset{(2)}{=}[\omega
]_{m_{\varepsilon }}\overset{not}{=}\beta _{m_{\varepsilon }}\in \Lambda
_{m_{\varepsilon }}(I)$, we get%
\begin{equation*}
d(\pi (\omega _{n}),\pi (\omega ))\leq diam(\overline{f_{\beta
_{m_{\varepsilon }}}(B)})=diam(f_{\beta _{m_{\varepsilon }}}(B))\overset{(3)}%
{<}\varepsilon \text{.}
\end{equation*}%
The last relation assures us that sequence $(\pi (\omega _{n}))_{n\in 
\mathbb{N}^{\ast }}$ converges to $\pi (\omega )$ and this shows that $\pi $
is continuous. $\square $

\bigskip

\textbf{Proposition 3.5.} \textit{Each diameter diminishing to zero IFS }$%
\mathcal{S}=((X,d),(f_{i})_{i\in I})$\textit{\ has attractor and }$A_{%
\mathcal{S}}=\pi (\Lambda (I))=\{a_{\omega }\mid \omega \in \Lambda (I)\}$.

\textit{Proof. }Let us note that:

i) 
\begin{equation*}
\pi (\Lambda (I))\overset{\text{Proposition 3.4 \& Remark 2.4, c)}}{\in }%
P_{cp}(X)\subseteq P_{b,cl}(X);
\end{equation*}

ii)%
\begin{equation*}
F_{\mathcal{S}}(\pi (\Lambda (I)))\overset{\text{Definition 2.10}}{=}%
\underset{i\in I}{\cup }f_{i}(\pi (\Lambda (I)))\overset{\text{Proposition
3.3}}{=}
\end{equation*}%
\begin{equation*}
=\underset{i\in I}{\cup }\pi (\tau _{i}(\Lambda (I)))=\pi (\underset{i\in I}{%
\cup }\tau _{i}(\Lambda (I)))=\pi (\Lambda (I))\text{.}
\end{equation*}

Moreover,%
\begin{equation*}
\underset{n\rightarrow \infty }{\lim }H_{d}(F_{\mathcal{S}}^{[n]}(B),\pi
(\Lambda (I)))=0\text{,}
\end{equation*}%
for each $B\in P_{b,cl}(X)$.

Indeed,\textit{\ }let us consider a fixed, but arbitrarily chosen $B\in
P_{b,cl}(X)$.\textit{\ }We have%
\begin{equation*}
H_{d}(F_{\mathcal{S}}^{[n]}(B),\pi (\Lambda (I)))=H_{d}(\underset{\omega \in
\Lambda (I)}{\cup }f_{[\omega ]_{n}}(B),\underset{\omega \in \Lambda (I)}{%
\cup }\{a_{\omega }\})\overset{\text{Proposition 2.2}}{\leq }
\end{equation*}%
\begin{equation*}
\leq \underset{\omega \in \Lambda (I)}{\sup }H_{d}(f_{[\omega
]_{n}}(B),\{a_{\omega }\})=\underset{\omega \in \Lambda (I)}{\sup }\max \{%
\underset{x\in B}{\sup }d(f_{[\omega ]_{n}}(x),a_{\omega }),\underset{x\in B}%
{\inf }d(f_{[\omega ]_{n}}(x),a_{\omega })\}=
\end{equation*}%
\begin{equation*}
\overset{a_{\omega }\in \overline{f_{[\omega ]_{n}}(B)}}{=}\underset{\omega
\in \Lambda (I)}{\sup }\underset{x\in B}{\sup }d(f_{[\omega
]_{n}}(x),a_{\omega })\text{,}
\end{equation*}%
for every $n\in \mathbb{N}^{\ast }$ and taking into account Proposition 3.2,
we conclude that $\underset{n\rightarrow \infty }{\lim }H_{d}(F_{\mathcal{S}%
}^{[n]}(B),\pi (\Lambda (I)))=0$. Consequently $\pi (\Lambda (I))$ is the
attractor of $\mathcal{S}$. $\square $

\bigskip

\textbf{Proposition 3.6.} \textit{For} \textit{each diameter diminishing to
zero IFS }$\mathcal{S}=((X,d),(f_{i})_{i\in I})$\textit{, we have }%
\begin{equation*}
\{a_{\overset{.}{\omega }}\}=Fix(f_{\omega })\text{,}
\end{equation*}%
\textit{for each} $\omega \in \Lambda ^{\ast }(I)\smallsetminus \{\lambda \}$%
.

\textit{Proof. }Since 
\begin{equation*}
f_{\omega }(a_{\overset{.}{\omega }})=f_{\omega }(\pi (\overset{.}{\omega }))%
\overset{\text{Proposition 3.3}}{=}\pi (\omega \overset{.}{\omega })=\pi (%
\overset{.}{\omega })=a_{\overset{.}{\omega }}\text{,}
\end{equation*}%
we conclude that 
\begin{equation}
\{a_{\overset{.}{\omega }}\}\subseteq Fix(f_{\omega })\text{,}  \tag{4}
\end{equation}%
for every $\omega \in \Lambda ^{\ast }(I)\smallsetminus \{\lambda \}$.

Now we prove that%
\begin{equation}
Fix(f_{\omega })\subseteq \{a_{\overset{.}{\omega }}\}\text{,}  \tag{5}
\end{equation}%
for every $\omega \in \Lambda ^{\ast }(I)\smallsetminus \{\lambda \}$.

Let us consider $z\in Fix(f_{\omega })$ (so $z\in \underset{n\in \mathbb{N}%
^{\ast }}{\cap }\{f_{[\overset{.}{\omega }]_{n\left\vert \omega \right\vert
}}(z)\}$). Then, in view of Definition 2.17, there exists $M_{\{z\}}\in
P_{b,cl}(X)$ such that $z\in M_{\{z\}}$ and $F_{\mathcal{S}%
}(M_{\{z\}})\subseteq M_{\{z\}}$. Hence 
\begin{equation*}
f_{i}(M_{\{z\}})\subseteq F_{\mathcal{S}}(M_{\{z\}})\subseteq M_{\{z\}}\text{%
,}
\end{equation*}%
for every $i\in I$, so $(\overline{f_{[\overset{.}{\omega }]_{n\left\vert
\omega \right\vert }}(M_{\{z\}})})_{n\in \mathbb{N}^{\ast }}$ is a
decreasing sequence and therefore 
\begin{equation*}
\underset{n\in \mathbb{N}^{\ast }}{\cap }\overline{f_{[\overset{.}{\omega }%
]_{n\left\vert \omega \right\vert }}(M_{\{z\}})}=\underset{n\in \mathbb{N}%
^{\ast }}{\cap }\overline{f_{[\overset{.}{\omega }]_{n}}(M_{\{z\}})}\overset{%
\text{Proposition 3.1}}{=}\{a_{\overset{.}{\omega }}\}\text{.}
\end{equation*}%
Consequently 
\begin{equation*}
z\in \underset{n\in \mathbb{N}}{\cap }\{f_{[\overset{.}{\omega }%
]_{n\left\vert \omega \right\vert }}(z)\}\subseteq \underset{n\in \mathbb{N}}%
{\cap }f_{[\overset{.}{\omega }]_{n\left\vert \omega \right\vert
}}(M_{\{z\}})\subseteq \underset{n\in \mathbb{N}}{\cap }\overline{f_{[%
\overset{.}{\omega }]_{n\left\vert \omega \right\vert }}(M_{\{z\}})}=\{a_{%
\overset{.}{\omega }}\}
\end{equation*}%
and $(4)$ is justified.

The relations $(4)$ and $(5)$ assure us that $Fix(f_{\omega })=\{a_{\overset{%
.}{\omega }}\}$. $\square $

\bigskip

\textbf{Proposition 3.7.} \textit{For} \textit{each diameter diminishing to
zero IFS }$\mathcal{S}=((X,d),(f_{i})_{i\in I})$\textit{, we have }%
\begin{equation*}
A_{\mathcal{S}}=\overline{\{a_{\overset{.}{\omega }}\mid \omega \in \Lambda
^{\ast }(I)\smallsetminus \{\lambda \}\}}\text{.}
\end{equation*}

\textit{Proof. }As $\{a_{\overset{.}{\omega }}\mid \omega \in \Lambda ^{\ast
}(I)\smallsetminus \{\lambda \}\}\subseteq A_{\mathcal{S}}$ and $A_{\mathcal{%
S}}$ is compact, it suffices to prove that 
\begin{equation*}
A_{\mathcal{S}}\subseteq \overline{\{a_{\overset{.}{\omega }}\mid \omega \in
\Lambda ^{\ast }(I)\smallsetminus \{\lambda \}\}}\text{.}
\end{equation*}

To this aim, let us consider $x\in A_{\mathcal{S}}$. For an arbitrary
neighborhood $V$ of $x$, there exists an open subset $D$ of $X$ such that $%
x\in D\subseteq V$. As $A_{\mathcal{S}}=\{a_{\omega }\mid \omega \in \Lambda
(I)\}$, there exists $\omega _{x}\in \Lambda (I)$ such that $x=a_{\omega
_{x}}$. Hence 
\begin{equation*}
\underset{n\in \mathbb{N}^{\ast }}{\cap }f_{[\omega _{x}]_{n}}(A_{\mathcal{S}%
})\overset{\text{Proposition 3.1}}{=}\{a_{\omega _{x}}\}\subseteq D
\end{equation*}%
and since the sequence $(f_{[\omega _{x}]_{n}}(A_{\mathcal{S}}))_{n\in 
\mathbb{N}^{\ast }}\subseteq P_{cp}(X)$ is decreasing and $D$ is open, there
exists $m\in \mathbb{N}$ such that 
\begin{equation}
f_{[\omega _{x}]_{m}}(A_{\mathcal{S}})\subseteq D  \tag{6}
\end{equation}%
(see, for example, Corollary 3.1.5 from [7]). Therefore 
\begin{equation*}
a_{\overset{\cdot }{_{[\omega _{x}]_{m}}}}\overset{\text{Proposition 3.6}}{=}%
f_{[\omega _{x}]_{m}}(a_{\overset{\cdot }{_{[\omega _{x}]_{m}}}})\in
f_{[\omega _{x}]_{m}}(A_{\mathcal{S}})\overset{\text{(6)}}{\subseteq }%
D\subseteq V\text{,}
\end{equation*}%
so 
\begin{equation*}
V\cap \{a_{\overset{.}{\omega }}\mid \omega \in \Lambda ^{\ast
}(I)\smallsetminus \{\lambda \}\}\neq \emptyset \text{.}
\end{equation*}%
Hence $x\in \overline{\{a_{\overset{.}{\omega }}\mid \omega \in \Lambda
^{\ast }(I)\smallsetminus \{\lambda \}\}}$ and the proof is done. $\square $

\bigskip

\textbf{4. THE\ MAIN\ RESULT}

\bigskip

\textbf{Theorem 4.1}. \textit{For an iterated function system} $\mathcal{S=}%
((X,d),(f_{i})_{i\in I})$, \textit{the following statements are equivalent:}

\textbf{1}. \textit{There exists a comparison function }$\varphi $\textit{\
such that} $\mathcal{S}$ \textit{is hyperbolic }$\varphi $\textit{%
-contractive.}

\textbf{2}. $\mathcal{S}$ \textit{is hyperbolic locally uniformly point
fibred.}

\textbf{3}. $\mathcal{S}$ \textit{is hyperbolic uniformly point fibred.}

\textbf{4}. $\mathcal{S}$ \textit{is a hyperbolic diameter diminishing to
zero iterated function system.}

\textbf{5}. $\mathcal{S}$\textit{\ has hyperbolic attractor and there exists
a continuous surjection }$\pi :\Lambda (I)\rightarrow A_{\mathcal{S}}$\ 
\textit{such that }%
\begin{equation*}
f_{i}\circ \pi =\pi \circ \tau _{i}\text{,}
\end{equation*}%
\textit{\ for all }$i\in I$\textit{.}

\textit{Proof}.

The argument used for the justification of Remark 3.1 from [15] ensures the
validity of 1)$\Rightarrow $3).

For 3)$\Rightarrow $2) see Remark 2.19, c).

For 2)$\Rightarrow $1) see Theorem 3.1 from [15].

1)$\Rightarrow $4) First of all, note that for every $B\in P_{b,cl}(X)$, we
have 
\begin{equation*}
\underset{\omega \in \Lambda _{n}(I)}{\max }diam(f_{\omega }(B))\leq \varphi
^{\lbrack n]}(diam(B))
\end{equation*}%
for every $n\in \mathbb{N}^{\ast }$, so, taking into account Remark 2.7, we
obtain that 
\begin{equation}
\underset{n\rightarrow \infty }{\lim }\underset{\omega \in \Lambda _{n}(I)}{%
\max }diam(f_{\omega }(B))=0\text{.}  \tag{7}
\end{equation}

In addition, since there exists a unique $A_{\mathcal{S}}\in P_{b,cl}(X)$
such that $F_{\mathcal{S}}(A_{\mathcal{S}})=A_{\mathcal{S}}$ and $\underset{%
n\rightarrow \infty }{\lim }F_{\mathcal{S}}^{[n]}(B)=A_{\mathcal{S}}$ for
every $B\in P_{b,cl}(X)$ (see Theorem 2.5 from [5]), one can consider the set%
\begin{equation}
M_{B}=\overline{A_{\mathcal{S}}\cup \underset{n\in \mathbb{N}}{\cup }%
F_{S}^{[n]}(B)}\supseteq B\text{.}  \tag{8}
\end{equation}

Note that as $\underset{n\rightarrow \infty }{\lim }F_{S}^{[n]}(B)=A_{%
\mathcal{S}}$, there exists $n_{1}\in \mathbb{N}$ and $\varepsilon >0$ such
that $F_{S}^{[n]}(B)\subseteq E_{\varepsilon }(A_{\mathcal{S}})$ for each $%
n\in \mathbb{N}$, $n\geq n_{1}$. Hence $\underset{n\in \mathbb{N},n\geq n_{1}%
}{\cup }F_{S}^{[n]}(B)\in P_{b}(X)$, so, via Definition 2.9, iv), we infer
that $A_{\mathcal{S}}\cup \underset{n\in \mathbb{N}}{\cup }F_{S}^{[n]}(B)\in
P_{b}(X)$. Consequently%
\begin{equation}
M_{B}\in P_{b,cl}(X)\text{.}  \tag{9}
\end{equation}

Moreover, we have%
\begin{equation*}
F_{S}(M_{B})=F_{S}(\overline{A_{\mathcal{S}}\cup \underset{n\in \mathbb{N}}{%
\cup }F_{S}^{[n]}(B)})=\overline{\underset{i\in I}{\cup }f_{i}(\overline{A_{%
\mathcal{S}}\cup \underset{n\in \mathbb{N}}{\cup }F_{S}^{[n]}(B)})}\subseteq
\end{equation*}%
\begin{equation*}
\overset{f_{i}\text{ continuous}}{\subseteq }\overline{\underset{i\in I}{%
\cup }\overline{f_{i}(A_{\mathcal{S}}\cup \underset{n\in \mathbb{N}}{\cup }%
F_{S}^{[n]}(B))}}\subseteq \overline{\overline{\underset{i\in I}{\cup }%
f_{i}(A_{\mathcal{S}}\cup \underset{n\in \mathbb{N}}{\cup }F_{S}^{[n]}(B))}}=
\end{equation*}%
\begin{equation*}
=\overline{F_{S}(A_{\mathcal{S}}\cup \underset{n\in \mathbb{N}}{\cup }%
F_{S}^{[n]}(B))}\subseteq \overline{F_{S}(A_{\mathcal{S}})\cup \underset{%
n\in \mathbb{N}}{\cup }F_{S}(F_{S}^{[n]}(B)))}\subseteq
\end{equation*}%
\begin{equation}
\subseteq \overline{A_{\mathcal{S}}\cup \underset{n\in \mathbb{N}}{\cup }%
F_{S}^{[n+1]}(B))}\subseteq M_{B}\text{.}  \tag{10}
\end{equation}

Therefore, in view of $(7)$, $(8)$, $(9)$ and $(10)$, $\mathcal{S}$ is a
hyperbolic diameter diminishing to zero iterated function system.

4)$\Rightarrow $5) See the results from Section 3.

5)$\Rightarrow $2) Let $d_{1}$ be the distance on $X$ whose existence is
assured by Remarks 2.19, a).

First of all let us note that $A_{\mathcal{S}}\in P_{cp}(X)$ (since $A_{%
\mathcal{S}}=\pi (\Lambda (I))$, $\pi $ is continuous and $\Lambda (I)$
compact -see Remark 2.4, c)-).

\textbf{Claim 1}.%
\begin{equation*}
\underset{n\rightarrow \infty }{\lim }\ \underset{\omega \in \Lambda _{n}(I)}%
{\max }diam(f_{\omega }(A_{\mathcal{S}}))=0\text{.}
\end{equation*}

\textit{Justification of Claim 1}. See Lemma 1.6 from [11].

\textbf{Claim 2}. For every $\omega \in \Lambda (I)$, the set $\underset{%
n\in \mathbb{N}^{\ast }}{\cap }f_{[\omega ]_{n}}(A_{\mathcal{S}})$ consists
on only one element denoted by $b_{\omega }$.

\textit{Justification of Claim 2}. Just use Cantor's theorem and Claim 1.

Let us consider $\varepsilon >0$ fixed, but arbitrarily chosen.

Then, via Claim 1, there exists $n_{0}\in \mathbb{N}$ such that 
\begin{equation}
diam(f_{\alpha }(A_{\mathcal{S}}))<\frac{\varepsilon }{3}\text{,}  \tag{11}
\end{equation}%
for every $\alpha \in \Lambda _{n_{0}}(I)$.

For every $\alpha \in \Lambda _{n_{0}}(I)$ there exists $\delta _{\alpha }>0$
such that 
\begin{equation}
f_{\alpha }(E_{\delta _{\alpha }}(A_{\mathcal{S}}))\subseteq E_{\frac{%
\varepsilon }{3}}(f_{\alpha }(A_{\mathcal{S}}))\text{.}  \tag{12}
\end{equation}

Indeed, the continuity of $f_{\alpha }$ assures us that for every $x\in A_{%
\mathcal{S}}$ there exists $\delta _{\alpha ,x}>0$ such that $f_{\alpha
}(B(x,2\delta _{\alpha ,x}))\subseteq B(f_{\alpha }(x),\frac{\varepsilon }{3}%
)$. In view of the compactness of $A_{\mathcal{S}}$ there exist $p\in 
\mathbb{N}^{\ast }$ and $x_{1},...,x_{p}\in A_{\mathcal{S}}$ such that $A_{%
\mathcal{S}}\subseteq \overset{p}{\underset{j=1}{\cup }}B(x_{j},\delta
_{\alpha ,x_{j}})$. Let us consider $\delta _{\alpha }=\min \{\delta
_{\alpha ,x_{1}},...,\delta _{\alpha ,x_{p}}\}>0$. If $u\in E_{\delta
_{\alpha }}(A_{\mathcal{S}})$, there exists $v\in A_{\mathcal{S}}$ such that 
$d_{1}(u,v)<\delta _{\alpha }$. In addition, there exists $j_{v}\in
\{1,...,p\}$ having the property that $d_{1}(x_{j_{v}},v)<\delta _{\alpha
,x_{_{j_{v}}}}$. Consequently we have $d_{1}(u,x_{j_{v}})\leq
d_{1}(u,v)+d_{1}(v,x_{j_{v}})<\delta _{\alpha }+\delta _{\alpha
,x_{_{j_{v}}}}<2\delta _{\alpha ,x_{_{j_{v}}}}$, so $f_{\alpha }(u)\in
f_{\alpha }(B(x_{j_{v}},2\delta _{\alpha ,x_{_{j_{v}}}}))\subseteq
B(f_{\alpha }(x_{j_{v}}),\frac{\varepsilon }{3})$, i.e. $f_{\alpha }(u)\in
E_{\frac{\varepsilon }{3}}(f_{\alpha }(A_{\mathcal{S}}))$. As $u$ was
arbitrarily chosen in $A_{\mathcal{S}}$, the justification of $(12)$ is done.

Let us consider 
\begin{equation*}
\delta =\min \{\delta _{\alpha }\mid \alpha \in \Lambda _{n_{0}}(I)\}\text{.}
\end{equation*}%
Since $\underset{n\rightarrow \infty }{\lim }F_{\mathcal{S}}^{[n]}(\overline{%
B(x,\eta )})=A_{\mathcal{S}}$ for every $x\in X$ and $\eta >0$, there exists 
$m_{0}\in \mathbb{N}^{\ast }$ such that 
\begin{equation}
F_{S}^{[m]}(\overline{B(x,\eta )})\subseteq E_{\delta }(A_{\mathcal{S}})%
\text{,}  \tag{13}
\end{equation}%
for every $m\in \mathbb{N}^{\ast }$, $m\geq m_{0}$. For an arbitrary $\omega
\in \Lambda (I)$ let us consider $[\omega ]_{n_{0}}\overset{not}{=}\alpha $.
Then, for $m\in \mathbb{N}$, $m\geq m_{0}$, if $[\omega ]_{n_{0}+m}=\alpha
\beta $, where $\beta \in \Lambda _{m}(I)$, via Claim 2, we have 
\begin{equation*}
d_{1}(f_{[\omega ]_{n_{0}+m}}(y),b_{\omega })=d_{1}(f_{\alpha }(f_{\beta
}(y)),b_{\omega })\leq \underset{u\in F_{S}^{m}(\overline{B(x,\eta )})}{\sup 
}d_{1}(f_{\alpha }(u),b_{\omega })\leq 
\end{equation*}%
\begin{equation*}
\overset{(13)}{\leq }\underset{u\in E_{\delta }(A_{\mathcal{S}})}{\sup }%
d_{1}(f_{\alpha }(u),b_{\omega })\overset{(12)}{\leq }\underset{v\in E_{%
\frac{\varepsilon }{3}}(f_{\alpha }(A_{\mathcal{S}})))}{\sup }%
d_{1}(v,b_{\omega })\leq 
\end{equation*}%
\begin{equation*}
\leq diam(E_{\frac{\varepsilon }{3}}(f_{\alpha }(A_{\mathcal{S}})))\leq 
\frac{2\varepsilon }{3}+diam(f_{\alpha }(A_{\mathcal{S}}))\overset{(11)}{<}%
\frac{2\varepsilon }{3}+\frac{\varepsilon }{3}=\varepsilon \text{,}
\end{equation*}%
for every $y\in \overline{B(x,\eta )}$, so%
\begin{equation*}
\underset{\omega \in \Lambda (I)}{\sup }\underset{y\in \overline{B(x,\eta )}}%
{\sup }d_{1}(f_{[\omega ]_{n_{0}+m}}(y),b_{\omega })\leq \varepsilon \text{,}
\end{equation*}%
for every $m\in \mathbb{N}$, $m\geq m_{0}$. Hence%
\begin{equation*}
\underset{n\rightarrow \infty }{\lim }\underset{\omega \in \Lambda (I)}{\sup 
}\underset{y\in \overline{B(x,\eta )}}{\sup }d_{1}(f_{[\omega
]_{n}}(y),b_{\omega })=0\text{,}
\end{equation*}%
i.e. $\mathcal{S}$ is hyperbolic uniformly point fibred. $\square $

\bigskip

\textbf{5. FINAL\ REMARKS}

\bigskip

\textbf{A.} Let us recall the concept of topologically contractive iterated
function system that was introduced by A. Mihail (see [19]) and by A.
Tetenov\ (see [20] and [21]) under a different name, namely self-similar
topological structure satisfying condition (P). It is part of the last
decades trend to establish purely topological conditions on an iterated
function system in order to guarantee the existence of attractors.

\bigskip

\textbf{Definition 5.1.} A\textit{\ pair }$((X,\tau ),(f_{i})_{i\in I})$%
\textit{\ is called a topologically contractive iterated function system if:}

\textit{i) }$(X,\tau )$\textit{\ is a topological space;}

\textit{ii) }$I$\textit{\ is a finite set;}

\textit{iii) }$f_{i}:X\rightarrow X$\textit{\ is continuous for every }$i\in
I$\textit{;}

\textit{iv) for every }$K\in P_{cp}(X)$\textit{\ there exists }$C_{K}\in
P_{cp}(X)$\textit{\ such that }$K\subseteq C_{K}$\textit{\ and }$\underset{%
i\in I}{\cup }f_{i}(C_{K})\subseteq C_{K}$\textit{;}

\textit{v) for every }$C\in P_{cp}(X)$\textit{\ such that }$\underset{i\in I}%
{\cup }f_{i}(C)\subseteq C$\textit{\ and every }$\omega \in \Lambda (I)$%
\textit{, the set \ }$\underset{n\in \mathbb{N}^{\ast }}{\cap }f_{[\omega
]_{n}}(C)$\textit{\ is a singleton.}

\bigskip

Particular cases of the above mentioned concept were considered by A. Edalat
(see [6]) under de name of weakly hyperbolic iterated function systems and
by B. Kieninger (see [12]) under the name of point fibred iterated function
systems.

\bigskip

A comparison between conditions i) and ii) from Definition 2.17 and
condition iv) from Definition 5.1 and between the conclusion of Proposition
3.1 and condition v) from Definition 5.1\ shows that the concept of diameter
diminishing to zero iterated function systems is a counterpart in terms of
metric spaces of the one of topologically contractive iterated function
system. Note that it is dealing with closed and bounded (not necessarly
compact) sets.

\bigskip

\textbf{B.}\ Even though Propositions 3.3, 3.4 and 3.5 follow from
Proposition 3.2 and Theorem 3.1 from [15] we presented their proofs as they
are elementary, while Theorem 3.1 from [15] is complicated and nontrivial.

\bigskip

\textbf{C.} For an iterated function system $\mathcal{S}=(((X,d),(f_{i})_{i%
\in I}))$ we can consider the following conditions:

1$^{^{\prime }}$) There exists a comparison function $\varphi $ such that $%
\mathcal{S}$ is $\varphi $-contractive.

2$^{^{\prime }}$) $\mathcal{S}$ is locally uniformly point fibred.

3$^{^{\prime }}$) $\mathcal{S}$ is uniformly point fibred.

4$^{^{\prime }}$) $\mathcal{S}$ is diameter diminishing to zero iterated
function system.

5$^{^{\prime }}$) $\mathcal{S}$ has an attractor and there exists a
continuous function $\pi :\Lambda (I)\rightarrow A_{\mathcal{S}}$ such that $%
\pi \circ \tau _{i}=f_{i}\circ \pi $ for all $i\in I$.

Actually the proof of Theorem 4.1 ensures the validity of the following
implications: 1$^{^{\prime }}$)$\rightarrow $4$^{^{\prime }}$), 4$^{^{\prime
}}$)$\rightarrow $3$^{^{\prime }}$), 3$^{^{\prime }}$)$\rightarrow $2$%
^{^{\prime }}$), 4$^{^{\prime }}$)$\rightarrow $5$^{^{\prime }}$) and 5$%
^{^{\prime }}$)$\rightarrow $2$^{^{\prime }}$).

We raise the following question: is 2$^{^{\prime }}$)$\rightarrow $4$%
^{^{\prime }}$) valid?

If this is true, then we get the equivalence of 2$^{^{\prime }}$), 3$%
^{^{\prime }}$), 4$^{^{\prime }}$) and 5$^{^{\prime }}$).

\bigskip

\textbf{Acknowledgement}. \textit{We want to thank the referees whose
generous and valuable remarks brought improvements to the paper (specially
by adding Section 5) and enhanced clarity.}

\bigskip

\textbf{References}

\bigskip

[1] R. Atkins, M. Barnsley, A. Vince, D. Wilson, A characterization of
hyperbolic affine iterated function systems, Topology Proc., \textbf{36}
(2010) 189--211.

[2] T. Banakh, W. Kubi\'{s}, N. Novosad, M. Nowak, F. Strobin, Contractive
function systems, their attractor and metrization, Topol. Methods Nonlinear
Anal., \textbf{46} (2015), 1029--1066.

[3] M. Barnsley, Fractals Everywhere, Academic Press, Boston, MA, 1988.

[4] M. Barnsley, A. Vince, Real projective iterated function systems, J.
Geom. Anal., \textbf{22} (2012), 1137-1172.

[5] D. Dumitru, Attractors of infinite iterated function systems containing
contraction type functions, An. \c{S}tiin\c{t}. Univ. Al. I. Cuza Ia\c{s}i
Mat. (N.S.), \textbf{59} (2013), 281-298.

[6] A. Edalat, Power domains and iterated function systems, Inf. Comput., 
\textbf{124} (1996), 182-197.

[7] R. Engelking, General Topology (Revised and completed edition),
Heldermann Verlag, Berlin, 1989.

[8] A. Granas, J. Dugundji, Fixed Point Theory, Springer-Verlag, New York,
2003.

[9] J. Hutchinson, Fractals and self similarity, Indiana Univ. Math. J., 
\textbf{30} (1981), 713--747.

[10] J. Jachymski, I. J\'{o}\'{z}wik, Nonlinear contractive conditions: a
comparison and related problems, Banach Center Publ. , \textbf{77} (2007),
123-146.

[11] A. Kameyama, Distances on topological self-similar sets and the
kneading determinants, J. Math. Kyoto Univ., \textbf{40} (2000), 603-674.

[12] B. Kieninger, Iterated function systems on compact Hausdorff spaces,
Ph.D. diss., University of Augsburg, Shaker-Verlag, Aachen, 2002.

[13] R. Miculescu, A. Mihail, Alternative characterization of hyperbolic
infinite iterated function systems, J. Math. Anal. Appl., \textbf{407}
(2013), 56-68.

[14] R. Miculescu, A. Mihail, On a question of A. Kameyama concerning
self-similar metrics. J. Math. Anal. Appl., \textbf{422}, 265--271 (2015).

[15] R. Miculescu, A. Mihail, A sufficient condition for a finite family of
continuous functions to be transformed into $\psi $-contractions, Ann. Acad.
Sci. Fenn., Math., \textbf{41} (2016), 51-65.

[16] R. Miculescu, A. Mihail, Remetrization results for possibly infinite
self-similar systems, Topol. Methods Nonlinear Anal., \textbf{47} (2016),
333-345.

[17] R. Miculescu, A. Mihail, A generalization of Istr\u{a}\c{t}escu's fixed
point theorem for convex contractions, Fixed Point Theory, \textbf{18}
(2017), 689-702.

[18] R. Miculescu, A. Mihail, A generalization for a finite family of
functions of the converse of Browder's fixed point theorem, Bull. Braz.
Math. Soc. (N.S.), \textbf{49} (2018), 673-698.

[19] A. Mihail, A topological version od iterated function systems, An. \c{S}%
tiin\c{t}. Al. I. Cuza, Ia\c{s}i, (S.N.), Matematica, \textbf{58} (2012),
105-120.

[20] M. Samuel, A. Tetenov, On attractors of iterated function systems in
uniform spaces.,Sib. \`{E}lektron. Mat. Izv., \textbf{14} (2017), 151--155.

[21] A. Tetenov, Semigroups satisfying P-condition and topological
self-similar sets, Sib. \`{E}lektron. Mat. Izv., \textbf{7} (2010), 461--464.

[22] S. Urziceanu, Alternative characterizations of AGIFSs having attractor,
Fixed Point Theory, \textbf{20} (2019), 729-740.

[23] A. Vince, Mobius iterated function systems, Trans. Amer. Math. Soc., 
\textbf{365} (2013), 491-509.

\bigskip

Radu Miculescu

Faculty of Mathematics and Computer Science

Transilvania University of Bra\c{s}ov

Iuliu Maniu Street, nr. 50, 500091, Bra\c{s}ov, Romania

E-mail: radu.miculescu@unitbv.ro

\bigskip

Alexandru Mihail

Faculty of Mathematics and Computer Science

University of Bucharest

Academiei Street 14, 010014, Bucharest, Romania

E-mail: mihail\_alex@yahoo.com

\end{document}